# Validity of Borodin and Kostochka Conjecture for classes of graphs without a single, forbidden subgraph on 5 vertices

Medha Dhurandhar

**Abstract:** Problem of finding an optimal upper bound for the chromatic no. of a graph is still open and very hard. **Borodin & Kostochka Conjecture** is still open and if proved will improve Brook's bound on Chromatic no. of a graph. Here we prove **Borodin & Kostochka Conjecture** for (1) ($P_4 \cup K_1$)-free (2) $P_5$-free (3) Chair-free graphs and 4) graphs with dense neighbourhoods. Certain known results follow as Corollaries.

**Introduction:**

In [1], [2], [3], [4] chromatic bounds for graphs are considered especially in relation with $\omega$ and $\Delta$. Gyárfás [5] and Kim [6] show that the optimal $\chi$-binding function for the class of $P_4 \cup K_1$-free graphs has order $\omega^2/\log(\omega)$. If we forbid additional induced subgraphs, the order of the optimal $\chi$-binding function drops below $\omega^2/\log(\omega)$. In 1941, Brooks' theorem stated that for any connected undirected graph $G$ with maximum degree $\Delta$, the chromatic number of $G$ is at most $\Delta$ unless $G$ is a complete graph or an odd cycle, in which case the chromatic number is $\Delta + 1$ [5]. In 1977, **Borodin & Kostochka** [6] conjectured that if $\Delta(G) \geq 9$, then $\chi(G) \leq \max\{\omega, \Delta-1\}$. In 1999, Reed proved the conjecture for $\Delta \geq 10^{14}$ [7]. Also D. W. Cranston and L. Rabern [8] proved it for claw-free graphs. We further proved it for a certain class of graphs [9]. Gupta and Pradhan [10] proved it for $\{C_4, P_5\}$-free graphs. Here we prove **Borodin & Kostochka** conjecture for graphs which are either $P_4 \cup K_1$-free or $P_5$-free or Chair-free and for graphs with dense neighbourhoods.

**Notation:** For a graph G, V(G), E(G), $\Delta$, $\omega$, $\chi$ denote the vertex set, edge set, maximum degree, size of a maximum clique, chromatic number of G resply. For $u \in V(G)$, $N(u) = \{v \in V(G) / uv \in E(G)\}$, and $\overline{N(u)} = N(u) \cup (u)$. If $S \subseteq V$, then <S> denotes subgraph of G induced by S. If C is some coloring of G and if $u \in V(G)$ is colored m in C, then u is called a m-vertex. If N(u) has a unique r-vertex say A, then r is called a unique color of u, A is called a unique r-vertex of u and if N(u) has more than one r-vertex, then r is called a repeat color of u. Also if P is a path in G s.t. vertices on P are alternately colored say i and j, then P is called a i-j bicolour path. All graphs considered henceforth are simple and undirected. For terms which are not defined herein we refer to Bondy and Murty [10].

We use the following **Lemma** in the subsequent **Main Results**.

**Lemma:** If G ($\neq K_{|V(G)|}$) is a smallest graph with $\Delta \geq 9$ and $\chi > \max\{\Delta-1, \omega\}$, then
1. $\chi = \Delta > \omega$
2. $\chi(G-u) = \chi-1 \; \forall \; u \in V(G)$.

Proof: As $G \neq K_{|V(G)|}$, $\Delta \geq \omega$. Now $\Delta-1 \geq \omega$ (else $\omega = \Delta$ and $\chi > \omega$. But as $G \neq C_{2n+1}$, by Brook's Theorem $\chi \leq \Delta = \omega$). Hence $\chi = \Delta > \omega$. Let $u \in V(G)$. Then $G-u \neq K_{|V(G)|-1}$ (else $\chi = \omega$). **If $\Delta(G-u) \geq 9$**, then by minimality, $\chi(G-u) \leq \max\{\omega(G-u), \Delta(G-u)-1\}$. Clearly if $\omega(G-u) \leq \Delta(G-u)-1$, then $\chi(G-u) \leq \Delta(G-u)-1 \leq \Delta-1 = \chi-1$ and otherwise $\chi(G-u) \leq \omega(G-u) \leq \omega < \Delta = \chi$. Thus $\chi(G-u) = \chi-1$. **Also if $\Delta(G-u) < 9$**, then as $G-u \neq C_{2n+1}$ (else $\chi \leq 4 < \Delta$) or $K_{|V(G)|-1}$ (else $\chi = \omega < \Delta$), by Brook's Theorem $\chi(G-u) \leq \Delta(G-u) < 9 \leq \Delta$. Thus always $\chi(G-u) = \chi-1 = \Delta-1 \; \forall \; v \in V(G)$.

**Main Result 1:** If G is a graph with $\Delta \geq 9$ s.t. for $u \in V(G)$, every $v \in N(u)$ is non-adjacent to at most three vertices in N(u), then $\chi \leq \max\{\Delta-1, \omega\}$.

Proof: Let if possible G be a smallest, connected graph with $\Delta \geq 9$ satisfying above condition and $\chi > \max\{\Delta-1, \omega\}$. Then by **Lemma**, $\chi = \Delta > \omega$, $\chi(G-u) = \chi-1 \; \forall \; u \in V(G)$. Let $u \in V(G)$ be s.t. deg $u = \Delta$. Let $S = \{1,..., \Delta-1\}$ be a $(\Delta-1)$-coloring of G-u. Then N(u) has at least $\Delta-2$ unique i-vertices say $A_i$ colored i ($1 \leq i \leq \Delta-2$) and at most a pair of vertices say X, Y with the same color $\Delta-1$. Clearly $A_i$ has a j-vertex for $1 \leq i \neq j \leq \Delta-2$ (else color $A_i$ by j, u by i) and $A_i$ has at most two j-vertices for any j, $1 \leq j \leq \Delta-2$ (else $A_i$ has a color $\alpha$ missing in $N(A_i)$, color $A_i$ by $\alpha$ and u by i).

**(a) If $A_iA_j \notin E(G)$, then for no m, $A_m$ is the only m-vertex of both $A_i$ and $A_j$ $1 \leq i, j, m \leq \Delta-2$.**
Let if possible $A_i$, $A_j$ both have $A_m$ as the only m-vertex. Then as $A_m$ has at the most one repeat color, w.l.g. let $A_i$ be the only i-vertex of $A_m$. Color $A_i$, $A_j$ by m, $A_m$ by i, u by j, a contradiction.

**(b) If $A_iA_j \notin E(G)$, then $A_i$, $A_j$ do not have more than two common adjacent vertices in N(u).**
**First** let $A_i$, $A_j$ be both adjacent to three unique vertices of u say $A_k$, $A_l$, $A_m$. As each of $A_i$, $A_j$ has at the most one repeat color, w.l.g. let $A_m$ be the only m-vertex of both $A_i$ and $A_j$, a contradiction to **(a)**. Next let $A_i$, $A_j$ be both adjacent to say $A_k$, $A_l$, X, $1 \leq i, j, k, l \leq \Delta-2$ where X has a repeat color. Then X is the unique ($\Delta-1$)-vertex of both $A_i$, $A_j$. If X has any color r missing in N(X), then color X by r, $A_i$ by $\Delta-1$, u by i, a contradiction. Hence X has no color missing and has either i or j as a unique colour in their neighbourhood. W.l.g. let i be the unique color for X. Then color $A_i$, $A_j$ by $\Delta-1$, X by i and u by j, a contradiction. Lastly let $A_i$, $A_j$ be both adjacent to say $A_k$, Y, X, $1 \leq i, j, k \leq \Delta-2$. As $A_i$, $A_j$ have two ($\Delta-1$)-vertices $A_k$ is their unique k-vertex, a contradiction to **(a)**.

**Case 1:** $\exists$ a ($\Delta-1$)-coloring of G-u s.t. $A_iA_j \notin E(G)$ where $A_i$, $A_j$ are unique vertices of u.
As $\Delta \geq 9$, $A_i$, $A_j$ have at least three common adjacent vertices in N(u), a contradiction to **(b).**

**Case 2:** In every ($\Delta-1$)-coloring of G-u, all vertices with unique colors in N(u) are adjacent.
Clearly $\Delta-1 \leq \omega \Rightarrow \Delta-1 = \omega \Rightarrow <u \cup \bigcup_{i=1}^{\Delta-2} A_i>$ is a maximum clique in G and $\{X, Y\} = N(u) - \bigcup_{i=1}^{\Delta-2} A_i$.
Now as every $v \in N(u)$ is non-adjacent to at most three vertices in N(u) and $\Delta \geq 9$, X and Y have at least three common adjacent vertices in N(u) say $A_i$ for $1 \leq i \leq 3$. W.l.g. let X be adjacent to more no. of $A_i$s than Y.

**Case 2.1:** X has no color missing in N(X).
Let $XA_k \notin E(G)$ for some $k > 3$. Now X has at the most one repeat color (else X has a color missing in N(X)), hence X and $A_k$ have a common unique vertex $A_1$. If X is the only ($\Delta-1$)-vertex of $A_1$, then color X, $A_k$ by 1, $A_1$ by $\Delta-1$, u by k, a contradiction. Hence let $A_1$ have two ($\Delta-1$)-vertices and a unique k-vertex. Color X, $A_k$ by 1 and $A_1$ by k. Then Y is the only ($\Delta-1$)-vertex of u, and $YA_j \in E(G)$ $\forall j \neq k$ (else we get **Case 1** with Y and $A_j$) and $YA_k \notin E(G)$. By assumption $XA_j \in E(G)$ $\forall j \neq k$. Again as $\Delta \geq 9$, $\exists A_p$ s.t. $A_p$ is the only p-vertex of X, Y and $A_k$. Color $A_p$ by k; X, Y, $A_k$ by p; u by ($\Delta-1$), a contradiction.

**Case 2.2:** Some color r is missing in N(X).
Then $XA_r \notin E(G)$ and $YA_j \in E(G)$ $\forall i \neq r$ (else color X by r and we get **Case 1** with Y and $A_i$) $\Rightarrow YA_r \notin E(G)$. Then by assumption $XA_j \in E(G)$ $\forall i \neq r$ and r is the only color missing in N(X). Again Y has no color t missing in N(Y) (else color X by r, Y by t and u by ($\Delta-1$)). Hence as $\Delta \geq 9$, $\exists A_p$ s.t. $A_p$ is the only p-vertex of X, Y and $A_r$. Color $A_p$ by ($\Delta-1$); X, Y, $A_r$ by p, u by r, a contradiction.

This proves **Case 2** and completes the proof of the **Main Result 1**.

**Main Result 2:** Let G be $P_4 \cup K_1$-free and $\Delta \geq 9$, then $\chi \leq \max\{\Delta-1, \omega\}$.
Proof: Let if possible G be a smallest, connected, $P_4 \cup K_1$-free graph with $\Delta \geq 9$ and $\chi > \max\{\Delta-1, \omega\}$. Then by **Lemma**, $\chi = \Delta > \omega$, $\chi(G-u) = \chi-1$ $\forall u \in V(G)$. Let $u \in V(G)$ be s.t. deg u = $\Delta$. Let $S = \{1,..., \Delta-1\}$ be a ($\Delta-1$)-coloring of G-u. Then N(u) has $\Delta-2$ unique vertices say $A_i$ colored i ($1 \leq i \leq \Delta-2$) and a pair of vertices say X, Y with the same color $\Delta-1$. Clearly $A_i$ has a j-vertex for $1 \leq i \neq j \leq \Delta-2$ (else color $A_i$ by j, u by i) and $A_i$ has at most two j-vertices for any j, $1 \leq j \leq \Delta-2$. Also as G is $P_4 \cup K_1$-free, every odd cycle of G is $C_5$.

**(a) If $A_iA_j \notin E(G)$, then for no m, $A_m$ is the only m-vertex of both $A_i$ and $A_j$ for $1 \leq i, j, m \leq \Delta-2$.**
Let if possible $A_i$, $A_j$ both have $A_m$ as the only m-vertex. Then as $A_m$ has at the most one repeat color, w.l.g. let $A_i$ be the only i-vertex of $A_m$. Color $A_i$, $A_j$ by m, $A_m$ by i, u by j, a contradiction.

**(b) If $A_iA_j \notin E(G)$, then $A_i$, $A_j$ do not have more than two common adjacent vertices in N(u).**

**First** let $A_i$, $A_j$ be both adjacent to say $A_k$, $A_l$, $A_m$ $1\le i, j, k, l, m\le\Delta-2$. As each of $A_i$, $A_j$ has at the most one repeat color, w.l.g. let $A_m$ be the only m-vertex of both $A_i$ and $A_j$, a contradiction to **(a)**. Next let $A_i$, $A_j$ be both adjacent to say $A_k$, $A_l$, X, $1\le i, j, k, l\le\Delta-2$. Then $A_i$, $A_j$ both have a unique $(\Delta-1)$-vertex X. If X has any color r missing in N(X), then color X by r, $A_i$ by $(\Delta-1)$, u by i, a contradiction. Hence X has no color missing and has either a unique i-vertex or j-vertex. W.l.g. let X have a unique i-vertex. Then color $A_i$, $A_j$ by $\Delta-1$, X by i, u by j, a contradiction. Lastly let $A_i$, $A_j$ be both adjacent to say $A_k$, Y, X, $1\le i, j, k\le\Delta-2$. As $A_i$, $A_j$ have two $(\Delta-1)$-vertices each $A_k$ is their only k-vertex, a contradiction to **(a)**.

**(c)** $A_i$ is non-adjacent to at the most two $A_k$s for $1\le i, k\le\Delta-2$.

Let if possible say $A_1A_k \notin E(G)$, $2\le k\le 4$. Now $A_k$ is adjacent to every 1-vertex ($\ne A_1$) of V(G) for $2\le k\le 4$ (else let n, n' ($\ne A_1$) be 1-vertices of G s.t. $nA_k \notin E(G)$ and $n'A_k \in E(G)$ $\Rightarrow$ <$A_1$, u, $A_k$, n', n> = $P_4\cup K_1$). Let $\{A_1, m, n, A_k\}$ be a 1-k path. Then n has two j-vertices for $2\le j\le 4$ (else if n has say one j-vertex, then $\exists$ $\{A_j, n', m', A_1\}$ j-1 path and <u, $A_1$, m', n', n> = $P_4\cup K_1$). Clearly this is true for all 1-vertices of G $\Rightarrow$ every 1-vertex has a color missing. Color each 1-vertex by its missing color and u by 1, a contradiction.

**Case 1:** $\exists$ a $(\Delta-1)$-coloring of G-u s.t. $A_iA_j \notin E(G)$ for some i, j $\in \{1,..,\Delta-2\}$.
As $\Delta \ge 9$, by **(c)** $\exists$ at least three vertices $A_k$s $(1\le i, k\le\Delta-2)$ adjacent to both $A_i$, $A_j$, a contradiction to **(b)**.

**Case 2:** In every $(\Delta-1)$-coloring of G-u, all vertices with unique colors in N(u) are adjacent.
Clearly $\Delta-1 \le \omega \Rightarrow \Delta-1 = \omega \ge 8 \Rightarrow <\bigcup_{i=1}^{\Delta-2} A_i \cup u>$ is a maximum clique in G and $\{X, Y\} = N(u) - \bigcup_{i=1}^{\Delta-2} A_i$.

**Claim: At most two vertices in $\bigcup_{i=1}^{\Delta-2} A_i$ are non-adjacent to both X and Y.**

Let if possible $A_i$ be non-adjacent to both X and Y for $1\le i\le 3$. Then if Z $\notin \{X, Y\}$ is any $(\Delta-1)$-vertex of G, $ZA_i \in E(G)$ for $1\le i\le 3$ (else if $ZA_i \notin E(G)$ and Z' is a $(\Delta-1)$-vertex of $A_i$, then <Z', $A_i$, u, X, Z> = $P_4\cup K_1$). Thus G has at most two $(\Delta-1)$-vertices $\notin \{X, Y\}$. Again at least one $(\Delta-1)$-vertex of $A_i$ has two i-vertices for $1\le i\le 3$. Hence clearly $\exists$ two $(\Delta-1)$-vertices Z, W $\notin \{X, Y\}$. W.l.g. say Z has two 1-vertices, 2-vertices and W has two 3-vertices. Let P = $\{A_1, Z, a_{11}, X\}$ be a 1-$(\Delta-1)$ bicolour path. Now $Ya_{11}$, $Wa_{11} \in E(G)$ (else <u, X, $a_{11}$, Z, Y/W> = $P_4\cup K_1$). Thus $a_{11}$ has four $(\Delta-1)$-vertices and has a color say r missing in N($a_{11}$). Color $a_{11}$ by r. Then X or Y has another 1-vertex say $a_{12}$. But as before it can be seen that $a_{12}$ has four $(\Delta-1)$-vertices and has a color say s missing in N($a_{12}$). Thus color X and Y by 1, u by $\Delta-1$, a contradiction.

Hence the **Claim** holds.

As $\Delta\ge 9$, w.l.g. let $A_i$ be adjacent to either X or Y for $1\le i\le 5$. Again w.l.g. let X be adjacent to more no. of $A_i$s than Y and hence let $XA_j \in E(G)$ for $1\le j\le 3$.

**Case 2.1:** X has no color missing in N(X).
Let $XA_k \notin E(G)$ for some k > 3. Now X has at the most one repeat color (else X has a color missing in N(X)) and clearly X and $A_k$ have a common unique vertex $A_n$ for $1\le n\le 3$. If X is the only $(\Delta-1)$-vertex of $A_n$, then color X, $A_k$ by n, $A_n$ by $\Delta-1$, u by k, a contradiction. Hence let $A_n$ have two $(\Delta-1)$-vertices and a unique k-vertex. Color X, $A_k$ by n and $A_n$ by k. Then Y is the only $(\Delta-1)$-vertex of u, $YA_j \in E(G)$ $\forall$ j $\ne$ k (else we get **Case 1** with Y and $A_j$) and $YA_k \notin E(G)$. By assumption $XA_j \in E(G)$ $\forall$ j $\ne$ k. Again as $\Delta\ge 9$, $\exists$ $A_p$ s.t. $A_p$ is the only p-vertex of X, Y and $A_k$. Color $A_p$ by k; X, Y and $A_k$ by p, u by $(\Delta-1)$, a contradiction.

**Case 2.2:** Some color r is missing in N(X).
Then $XA_r \notin E(G)$ and $YA_j \in E(G)$ $\forall$ i $\ne$ r (else color X by r and we get **Case 1** with Y and $A_i$) $\Rightarrow YA_r \notin E(G)$. Again Y has no color t missing in N(Y) (else color Y by t and u by $(\Delta-1)$). Hence as $\Delta\ge 9$, $\exists$ $A_p$ s.t. $A_p$ is the only p-vertex of Y and $A_r$. Now $A_p$ has either a unique $(\Delta-1)$-vertex or r-vertex. Color $A_p$ by $(\Delta-1)$/r; Y and $A_r$ by p, u by r/$(\Delta-1)$, a contradiction.

This proves **Case 2** and completes the proof of the **Main Result 2**.

Our results in [9] follow as corollaries to the **Main Result 2**.

**Corollary 1: Borodin & Kostochka** conjecture is true for $P_3 \cup K_1$-free graphs.
**Corollary 2: Borodin & Kostochka** conjecture is true for $K_2 \cup \overline{K_2}$-free graphs.
**Corollary 3: Borodin & Kostochka** conjecture is true for $3K_1$-free graphs.

**Main Result 3:** Let G be Chair-free and $\Delta \geq 9$, then $\chi \leq \max\{\Delta-1, \omega\}$.
Proof: Let if possible G be a smallest, connected, Chair-free graph with $\Delta \geq 9$ and $\chi > \max\{\Delta-1, \omega\}$. Then by **Lemma**, $\chi = \Delta > \omega$, $\chi(G-u) = \chi-1 \ \forall \ u \in V(G)$. Let $u \in V(G)$ be s.t. deg $u = \Delta$. Let $S = \{1,..., \Delta-1\}$ be a $(\Delta-1)$-coloring of G-u. Then N(u) has at least $\Delta-2$ vertices say $A_i$ with unique colors i ($1 \leq i \leq \Delta-2$) and at most a pair of vertices say X, Y with the same color $\Delta-1$. Clearly $A_i$ has a j-vertex for $1 \leq i \neq j \leq \Delta-2$ (else color $A_i$ by j, u by i) and $A_i$ has at most two j-vertices for any j, $1 \leq j \leq \Delta-2$.

**(a) For no m, $A_m$ is the only m-vertex of both $A_i$ and $A_j$ for $1 \leq i, j, m \leq \Delta-2$.**
Let if possible $A_i$, $A_k$ both have $A_m$ as the only m-vertex. Then as $A_m$ has at the most one repeat color, w.l.g. let $A_i$ be the only i-vertex of $A_m$. Color $A_i$, $A_k$ by m, $A_m$ by i, u by k, a contradiction.

**(b) If $A_iA_j \notin E(G)$, then $A_i$ ($A_j$) has a unique j-vertex (i-vertex).**
If Z, W are the two j-vertices of $A_i$, then $<A_i, Z, W, u, A_j>$ = Chair, a contradiction.

**(c) $A_i$, $A_j$ do not have more than two common adjacent vertices in N(u).**
**First** let $A_i$, $A_j$ be both adjacent to say $A_k$, $A_l$, $A_m$ $1 \leq i, j, k, l, m \leq \Delta-2$. As each of $A_i$, $A_j$ has at the most one repeat color, w.l.g. let $A_m$ be the only m-vertex of both $A_i$ and $A_j$. This is a contradiction to **(a)**. Next let $A_i$, $A_j$ be both adjacent to say $A_k$, $A_l$, X, $1 \leq i, j, k, l \leq \Delta-2$. Then by **(a)**, $A_i$, $A_j$ both have a unique x-vertex X. If X has any color r missing in N(X), then color X by r, $A_i$ by x, u by i, a contradiction. Hence X has no color missing and has either a unique i-vertex or j-vertex. W.l.g. let X have a unique i-vertex. Then color $A_i$, $A_j$ by $\Delta-1$, X by i, u by j, a contradiction. Lastly let $A_i$, $A_j$ be both adjacent to say $A_k$, Y, X, $1 \leq i, j, k \leq \Delta-2$. As $A_i$, $A_j$ have two $(\Delta-1)$-vertices each they have a unique k-vertex, a contradiction to **(a)**.

**(d) $A_i$ is non-adjacent to at the most two $A_k$s for $1 \leq i, k \leq \Delta-2$.**
Let if possible $A_1A_k \notin E(G)$ for $2 \leq k \leq 4$. By **(b)**, each $A_k$ has a unique i-vertex and that i-vertex has two k-vertices $\Rightarrow$ at least two $A_k$'s have different i-vertices. W.l.g. let $ZA_2, WA_3 \in E(G)$ where Z, W are i-vertices. Again by **(b)**, w.l.g. let $ZA_4 \notin E(G)$. Now $A_2A_3, A_2A_4 \in E(G)$ (else $<u, A_1, A_2, A_3, Z>$ = Chair). As $A_2$ has at most one repeat color, w.l.g. let $A_2$ have a unique 3-vertex. Consider a 1-3 bi-color component R containing $A_1, A_3$. Let P be a $A_1$-$A_3$ path in R.

**Claim:** R = P.
Let if possible V be the **First** vertex on P with degree more than two in R. By **(b)**, $V \neq A_1$. Let M, N be the vertices preceding V on P and L be the vertex following V on P. Let S be the third vertex having same color as N and L s.t. L is adjacent to V. But then either $<V, M, L, S, N>$ = Chair or M has degree at least three in R, a contradiction. Hence the **Claim** holds.

Alter colors along P, color $A_2$ by 3, u by 2, a contradiction.

**Case 1:** $\exists$ a $(\Delta-1)$-coloring of G-u s.t. $A_iA_j \notin E(G)$ for some i, j $\in \{1,.., \Delta-2\}$.
As $\Delta \geq 9$, by **(c)** $\exists$ at least three vertices $A_k$s ($1 \leq i, k \leq \Delta-2$) adjacent to both $A_i, A_j$, a contradiction to **(b).**

**Case 2:** In every $(\Delta-1)$-coloring of G-u, all vertices with unique colors in N(u) are adjacent.
Clearly $\Delta-1 \leq \omega \Rightarrow \Delta-1 = \omega \geq 8 \Rightarrow <\bigcup_{i=1}^{\Delta-2} A_i>$ is a maximum clique in G-u and $\{X, Y\} = N(u) - \bigcup_{i=1}^{\Delta-2} A_i$.

As G is Chair-free, clearly every $A_i$ is adjacent to X or Y. Let X have more $A_k$s adjacent than Y. As $\Delta \geq 9$, $\exists$ at least four $A_k$s adjacent to X. Let $XA_s \notin E(G) \Rightarrow YA_s \in E(G)$. Let $YA_t \notin E(G)$. If X has a color say s missing, then we get **Case 1** with Y and $A_t$, a contradiction. Hence X has no color missing and X and $A_s$ have at least two unique colors say i, j. If $A_i$ has a unique ($\Delta$-1)-vertex X, then color X and $A_s$ by i, $A_i$ by $\Delta$-1, u by i, a contradiction. If $A_i$ has two ($\Delta$-1)-vertices, then color X and $A_s$ by i, $A_i$ by s. Thus Y is a unique ($\Delta$-1)-vertex of N(u) and we get **Case 1** with Y and $A_t$, a contradiction.

This proves **Case 2** and completes the proof of the **Main Result 3**.

D. W. Cranston and L. Rabern [8] result for claw-free graphs follows as a Corollary.
**Corollary:** Borodin & Kostochka conjecture is true for claw-free graphs.

**Main Result 4:** Let G be $P_5$–free and $\Delta \geq 9$, then $\chi \leq \max \{\Delta-1, \omega\}$.
Proof: Let if possible G be a smallest, connected, $P_5$-free graph with $\Delta \geq 9$ and $\chi > \max \{\Delta-1, \omega\}$. Then by **Lemma**, $\chi = \Delta > \omega$, $\chi(G-u) = \chi-1 \ \forall \ u \in V(G)$. Let $u \in V(G)$ be s.t. deg u = $\Delta$. Let S = {1,..., $\Delta$-1, $\Delta$} be a $\Delta$-coloring of G with only u colored $\Delta$. Then N(u) has $\Delta$-2 vertices say $A_i$ with unique colors i ($1 \leq i \leq \Delta-2$) and a pair of vertices say X, Y with the same color $\Delta$-1. Clearly $A_i$ has a j-vertex for $1 \leq i \neq j \leq \Delta-2$ (else color $A_i$ by j, u by i) and $A_i$ has at most two j-vertices for any j, $1 \leq j \leq \Delta-2$. Also as G is $P_4 \cup K_1$-free, every odd cycle of G is $C_5$.

**(a)** For no m, $A_m$ is the only m-vertex of both $A_i$ and $A_j$ for $1 \leq i, j, m \leq \Delta-2$.
Let if possible $A_i$, $A_k$ both have $A_m$ as the only m-vertex. Then as $A_m$ has at the most one repeat color, w.l.g. let $A_i$ be the only i-vertex of $A_m$. Color $A_i$, $A_k$ by m, $A_m$ by i, u by k, a contradiction.

**(b)** $A_i$, $A_j$ do not have more than two common adjacent vertices in N(u).
**First** let $A_i$, $A_j$ be both adjacent to say $A_k$, $A_l$, $A_m$ $1 \leq i, j, k, l, m \leq \Delta-2$. As each of $A_i$, $A_j$ has at the most one repeat color, w.l.g. let $A_m$ be the only m-vertex of both $A_i$ and $A_j$. This is a contradiction to **(a)**. Next let $A_i$, $A_j$ be both adjacent to say $A_k$, $A_l$, X, $1 \leq i, j, k, l \leq \Delta-2$. Again as before $A_i$, $A_j$ both have a unique x-vertex X. If X has any color r missing in N(X), then color X by r, $A_i$ by x, u by i, a contradiction. Hence X has no color missing and has either a unique i-vertex or j-vertex. W.l.g. let X have a unique i-vertex. Then color $A_i$, $A_j$ by $\Delta$-1, X by i, u by j, a contradiction. Lastly let $A_i$, $A_j$ be both adjacent to say $A_k$, Y, X, $1 \leq i, j, k \leq \Delta-2$. As $A_i$, $A_j$ have two ($\Delta$-1)-vertices each they have a unique k-vertex, a contradiction to **(a)**.

**(c)** $A_i$ is non-adjacent to at the most two $A_k$s for $1 \leq i, k \leq \Delta-2$.
Let if possible $A_1 A_k \notin E(G)$ for $2 \leq k \leq 4$. Now $\exists$ a 1-vertex say n of $A_2$ s.t. n has no color missing in N(n) (else color each such n by the missing color, $A_2$ by 1, u by 2). Let {$A_1$, m, n, $A_2$} be a 1-2 path from $A_1$ to $A_2$. As n has two 2-vertices, w.l.g. let n have at most one 3-vertex (else if n has two 2-vertices and 3-vertices, then n has a color say r missing). Let {$A_1$, s, t, $A_3$} be a 1-3 path from $A_1$ to $A_3$. Now ns $\in$ E(G) iff $nA_3 \in$ E(G) (else <n, s/$A_3$, $A_1$, u, $A_3$/s> = $P_5$). Hence ns, $nA_3 \notin$ E(G). Let v be a 3-vertex of n. Again $sA_2 \in$ E(G) (else <n, $A_2$, u, $A_1$, s> = $P_5$) and $vA_2 \notin$ E(G) (else $A_2A_3$, $A_2t \notin$ E(G) and <n, $A_2$, u, $A_3$, t> = $P_5$) $\Rightarrow A_2A_3 \in$ E(G) (else <v, n, $A_2$, u, $A_3$> = $P_5$). Thus $A_2$ has two 3-vertices and hence n is its only 1-vertex. But then <v, n, $A_2$, $A_3$, t> = $P_5$, a contradiction.

**Case 1:** $\exists$ a ($\Delta$-1)-coloring of G-u s.t. $A_iA_j \notin$ E(G) for some i, j $\in$ {1,.., $\Delta$-2}.
As $\Delta \geq 9$, by **(c)** $\exists$ at least three vertices $A_k$s ($1 \leq i, k \leq \Delta-2$) adjacent to both $A_i$, $A_j$, a contradiction to **(b).**

**Case 2:** In every ($\Delta$-1)-coloring of G-u, all vertices with unique colors in N(u) are adjacent.
Clearly $\Delta-1 \leq \omega \Rightarrow \Delta-1 = \omega \geq 8 \Rightarrow < \bigcup_{i=1}^{\Delta-2} A_i >$ is a maximum clique in G-u and {X, Y} = N(u)- $\bigcup_{i=1}^{\Delta-2} A_i$.

**Claim: At most two vertices in $\bigcup_{i=1}^{\Delta-2} A_i$ are non-adjacent to both X and Y.**

Let if possible $A_i$ be non-adjacent to both X and Y for $1 \leq i \leq 3$. Then $\exists$ a i-($\Delta$-1) path {$A_i$, $a_{(\Delta-1)1}$, $a_{i1}$, X/Y}. Again $a_{i1}Y/X \in$ E(G) for $1 \leq i \leq 3$ (else <Y/X, u, $A_i$, $a_{(\Delta-1)1}$, $a_{i1}$> = $P_5$). **First let $A_i$ have two ($\Delta$-1)-vertices W and Z**. As G is $P_5$-free, $Wa_{i1}$, $Za_{i1} \in$ E(G). Thus $a_{i1}$ has four ($\Delta$-1)-vertices and hence a color say r missing in N($a_{i1}$). This is true for every i-vertex adjacent to X or Y. By coloring all such $a_{i1}$s by the missing colors, coloring X and Y by i and u by Δ-1, we get a contradiction. **Next let every $A_i$ have a unique (Δ-1)-vertex.** Clearly ∃ i, j (1≤i, j≤3) s.t. $A_i$, $A_j$ have different (Δ-1)-vertices (else the common (Δ-1)-vertex of all $A_i$s has a color say r missing. Color that vertex by r, $A_i$ by Δ-1, u by i) say Z, W resply. Consider a j-(Δ-1) path {$A_j$, W, $a_{j1}$, X} from $A_j$ to X. Then $Za_{j1}$ ∈ E(G) and $a_{j1}$ has a color missing in $N(a_{j1})$. As before, we get a contradiction.

Thus the **Claim** holds.

As Δ≥9, w.l.g. let $A_i$ be adjacent to either X or Y for 1≤i≤5. Again w.l.g. let X be adjacent to more no. of $A_i$s than Y and hence let $XA_j$ ∈ E(G) for 1≤j≤3.

**Case 2.1:** X has no color missing in N(X).
Let $XA_k$ ∉ E(G) for some k > 3. Now X has at the most one repeat color and clearly X and $A_k$ have a common unique vertex $A_n$ for 1≤n≤3. If X is the only (Δ-1)-vertex of $A_n$, then color X, $A_k$ by n, $A_n$ by Δ-1, u by k, a contradiction. Hence let $A_n$ have two (Δ-1)-vertices and a unique k-vertex. Color X, $A_k$ by n and $A_n$ by k. Then Y is the only (Δ-1)-vertex of u, $YA_j$ ∈ E(G) ∀ j ≠ k (else we get **Case 1** with Y and $A_j$) and $YA_k$ ∉ E(G). By assumption $XA_j$ ∈ E(G) ∀ j ≠ k. Again as Δ≥9, ∃ $A_p$ s.t. $A_p$ is the only p-vertex of X, Y and $A_k$. Color $A_p$ by k; X, Y and $A_k$ by p, u by (Δ-1), a contradiction.

**Case 2.2:** Some color r is missing in N(X).
Then $XA_r$ ∉ E(G) and $YA_j$ ∈ E(G) ∀ i ≠ r (else color X by r and we get **Case 1** with Y and $A_i$) ⇒ $YA_r$ ∉ E(G). Again Y has no color t missing in N(Y) (else color Y by t and u by (Δ-1)). Hence as Δ≥9, ∃ $A_p$ s.t. $A_p$ is the only p-vertex of Y and $A_r$. Now $A_p$ has either a unique (Δ-1)-vertex or r-vertex. Color $A_p$ by (Δ-1)/r; Y and $A_r$ by p, u by r/(Δ-1), a contradiction.

This proves **Case 2** and completes the proof of the **Main Result 4**.

Gupta and Pradhan [10] result for {$C_4$, $P_5$}-free graphs follows as a Corollary.
**Corollary:** Borodin & Kostochka conjecture is true for {$C_4$, $P_5$}-free graphs.